\documentclass[a4paper, 12pt]{article}
\usepackage{adjustbox}
\usepackage{graphicx}
\usepackage{amsfonts}
\usepackage{amsbsy}
\usepackage{amssymb}
\usepackage[fleqn]{amsmath}
\usepackage{multicol}
\usepackage{longtable,ltxtable,booktabs}
\usepackage{blkarray}
\usepackage{afterpage}
\usepackage{float}
\usepackage{caption}
\usepackage{subcaption}
\usepackage{multirow}
\usepackage{pdflscape}
\usepackage{latexsym}
\usepackage{rotating}
\usepackage{enumerate}

\usepackage{mathtools}

\usepackage{epstopdf}
\usepackage{caption}
\usepackage{xcolor}
\usepackage[nodisplayskipstretch]{setspace}
\usepackage{amsmath}
\usepackage{longtable}
\usepackage[margin=.8 in]{geometry}

\usepackage{verbatim}

\usepackage{mathrsfs}
\usepackage{amsthm}
\usepackage{setspace}
 \usepackage{pdflscape}
\usepackage{colortbl}
\usepackage{latexsym}
\usepackage{amsfonts}
\usepackage{epsfig}
\usepackage{hyperref}
\usepackage{graphics}
\usepackage{float}
\theoremstyle{definition}
\newtheorem{thm}{Theorem}[section]

\newtheorem{defn}{Definition}[section]
\newtheorem{lemma}{Lemma}[section]

\newtheorem{cor}{Corollary}[section]

\theoremstyle{definition}

\title{ GVZ-groups with two character degrees}

\date{{}}\author{Nabajit Talukdar\thanks{Corresponding author.Email address : ntalukdar2000@yahoo.co.in}\\
	\small and\\
	Kukil Kalpa Rajkhowa\thanks{ E-mail address : kukilrajkhowa@yahoo.com}\\
	\small Department of Mathematics\\
	\small Cotton University\\
	\small Guwahati-781001, India}

\begin{document}
\maketitle

\begin{abstract}
 We investigate the finite groups $G$ for which $\chi(1)^{2}=|G:Z(\chi)|$ for all characters $\chi \in Irr(G)$ and $|cd(G)|=2$, where $cd(G)=\{\chi(1)| \chi \in Irr(G)\}$. We obtain some alternate characterizations of these groups and we obtain some information regarding the structures of these groups.
 \end{abstract}

 $\mathbf{2010\;Mathematics\;Subject\;Classification:}\;20C15$
	
	$\mathbf{Keywords\;and\;Phrases:}$ GVZ-groups, Character degrees, $p$-groups

\section{Introduction}
In this paper, all groups are finite and we write $Irr(G)$ and $nl(G)$ for the set of irreducible  characters and the set of non-linear irreducible  characters of $G$ respectively. $c(G)$ denotes the nilpotency class of $G$. For the terminologies not defined here follow the Isaacs' book \cite{isaacs1994character} and Rotman's book \cite{rotman2012introduction}. An important question in group theory is what information can we get from the set of character degrees about the structure of the group. There are many research articles attempting to answer this question by examining the groups whose irreducible characters have some distinctive properties like "the group has only two character degrees", "the centres of the irreducible characters form a chain with respect to inclusion" etc. In \cite{nenciu2016nested}, A. Nenciu obtained various characterizations of the finite groups for which $\chi(1)^{2}=|G:Z(\chi)|$ for all characters $\chi \in Irr(G)$ and $\{Z(\chi):\chi \in Irr(G)\}$
 is a chain, where $Z(\chi)=\{g\in G: |\chi(g)|=\chi(1) \}$. In \cite{lewis2020groups} Mark L. Lewis obtained the characterizations of the groups for which $\{Z(\chi):\chi \in Irr(G)\}$
 is a chain with respect to inclusion and obtained some information regarding the structure of these groups. In this paper we examine the GVZ-groups having two character degrees. The concept of GVZ-group first appeared in \cite{murai1994} under the name "groups of Ono type". The term  GVZ-group was introduced by A. Nenciu.

 \begin{defn}
 \cite{nenciu2016nested}
A nonabelian group $G$ is called a generalized VZ-group (GVZ for short) if for all $\chi\in Irr(G)$ we have $\chi(g)=0$ for all $g\in G\setminus Z(\chi)$.
 \end{defn}
 By Corollary 2.30 \cite{isaacs1994character} we get that if $G$ is a GVZ-group then $\chi(1)^{2}=|G:Z(\chi)|$. If $G$ is a GVZ-group with two character degrees then it follows that $|Z(\chi)|$ is constant for all $\chi \in nl(G)$. In \cite{lewis2020groups} Mark L. Lewis obtained various results about the factor groups for a group where centres of the irreducible characters form a chain. In this paper, we examine the GVZ-groups having two character degrees. It is to be noted that if a GVZ-group has only two character degrees, then the centres of non-linear irreducible characters of the group are of equal size. We obtain the following results.
 \begin{thm}
 Let $G$ be a  GVZ-group with two character degrees. Then for any $\chi\in nl(G)$, the groups 
 $G/Z(\chi)$, $G/Z(G)$ and $G'/[Z(\chi),G]$
  are elementary Abelian and $exp\ G/Z(\chi)=exp\ G/Z(G)=exp\ G'/[Z(\chi),G]$.
 \end{thm}

 \begin{thm}
Let $G$ be a GVZ group with two character degrees. 
Let $N$  be a normal subgroup of $G$ such that $G'\not \subseteq N$. Then  $\chi \in nl(G/N)$  if and only if $N[Z(\chi),G]<NG'$.
 %Let $N$  be a normal subgroup of $G$ such that $G'\not \subseteq N$. Then  $\chi \in nlIrr(\frac{G}{N})$  iff $N[Z(\chi),G]<NG'$.
 \end{thm}

We also obtain a characterization of  GVZ-groups with two character degrees by the orders of the centres of some factor groups of the group.

\begin{thm}
Let $G$ be a   GVZ-group. Then the following conditions are equivalent:
  \begin{enumerate}[(i)]
   \item $|cd(G)|=2$.

   \item There is a collection $\{X_{1}, X_{2}, \ldots , X_{n}\}$ of normal subgroups of $G$ of equal order such that  every normal subgroup $N$ of $G$ with $G'\not \subseteq N$ and $[Z(\chi),G ]\subseteq N$ for some $\chi \in \ nl(G)$ satisfies $Z(G/N)=X_{i}/N$ for some $1\leq i \leq n$.

   \item For any two normal subgroups $N$ and $M$ of $G$ such that $G'\not \subseteq N,M$ and $[Z(\chi), G] \subseteq N$, $[Z(\phi), G] \subseteq M$ for some $\chi, \phi \in \ nl(G)$ if we define $Z_{N}/N=Z(G/N)$, $Z_{M}/M=Z(G/M)$ then $|Z_{N}|=|Z_{M}|$.

   Moreover the collection $\{X_{1}, X_{2}, \ldots , X_{n}\}$ of normal subgroups in (ii) is the collection of centres of non-linear irreducible characters of $G$.

  \end{enumerate}

  \end{thm}

  The characterizations of $p$-groups from the set of  character degrees
  are of special interest. For example the finite $p$-groups with two character degrees  where all the non-linear irreducible characters are faithful have been studied in \cite{di1992some} and \cite{doostie2012finite}. The finite $p$- groups $G$ having character degree set $\{1, |G:Z(G)|^{\frac{1}{2}}\}$ have been characterized in \cite{fernandez2001groups}. Generalizing this result we obtain the following theorem.
  \begin{thm}
Let $G$ be a non-Abelian $p$ group with $|cd(G)|=2$ .Then the following conditions are equivalent:
\begin{enumerate}[(i)]
    \item  $G$ is a GVZ group.
   
    %\item $Cl(x[Z(\chi), G])=xG'/[Z(\chi), G]$ in the group $G/[Z(\chi),G]$ for all $x \not \in  Z(\chi)$, $\chi \in \ nl(G)$.

     \item $G/[Z(\chi),G]$ is isoclinic to a semi extraspecial $p$-group for any $\chi \in nl(G)$.
    
    \item $Z(G/N)=Z(\chi)/N$ for any normal subgroup $N$ of $G$ such that $G'\not \leq N$ and $[Z(\chi),G]\leq N$ for some $\chi \in \ nl(G)$.
    
\end{enumerate}
  \end{thm}

\section{Preliminaries}
The aim of this section is to prove some essential results which will be applied in the proofs of our main results. We state the following lemma.

\begin{lemma}
\label{lemma1}
Let $G$ be a GVZ-group with two character degrees. Then the nilpotency class of $G$ is two.
\begin{proof}
From Theorem B in \cite{burkett2021gvz} we get that if $G$ is a GVZ-group then $c(G)\leq |cd(G)|$.Since $|cd(G)|=2$ and $G$ is a non-Abelian group it follows that $c(G)=2$. 
\end{proof}
\end{lemma}

In the following Lemma we get  an essential condition that the centres of two non-linear irreducible characters of a GVZ-group with two character degrees are equal.

\begin{lemma}
\label{lemma2}
Let $G$ be a GVZ group with two character degrees. Let $\phi, \chi\in nl(G)$ . If $[Z(\phi), G]\leq ker\chi$, then $Z(\phi)=Z(\chi)$.
\begin{proof}
Suppose $[Z(\phi), G]\leq ker\chi$.
Then $\frac{Z(\phi)ker\chi}{\ker\chi}\subseteq Z(\frac{G}{ker\chi})=\frac{Z(\chi)}{ker\chi}$.\\
This gives $Z(\phi)\subseteq Z(\chi)$ . Since $|Z(\chi)|$ is constant for all $\chi\in nl(G)$, we get $Z(\phi)=Z(\chi)$ .
\end{proof}
\end{lemma}

\begin{cor}
Let $G$ be a GVZ group with two character degrees. Let  $N\unlhd G$
and $ \chi\in nl(G)$ be such that $G'\not \subseteq N $ and $[Z(\chi), G]\leq N$. Then $N\leq Z(\chi)$.

\begin{proof}
Let $\phi\in nl(G)$ be such that $N\leq ker\phi$. Then $N\leq Z(\phi)$ and since $[Z(\chi), G]\leq ker\phi$, it follows that $Z(\chi)=Z(\phi)$. Thus we get that
$N\leq Z(\chi)$.

\end{proof}

\end{cor}

The following lemma gives a necessary and sufficient condition to find the centres of the non-linear irreducible characters of a GVZ-group with two character degrees.

\begin{lemma}
\label{lemma3}
Let $G$ be a GVZ group with two character degrees. 
 Let $\{X_{1}, X_{2}, \ldots , X_{n}\}$ be the collection of centres of non-linear irreducible characters. Let $\chi \in nl(G)$. Then $Z(\chi)=X_{i}$ if and only if $[X_{i}, G]\leq ker\chi$.

\begin{proof}
If $Z(\chi)=X_{i}$, then $[X_{i}, G]=[Z(\chi),G]\leq ker \chi$. Conversely if $[X_{i}, G]\leq ker\chi$ then $X_{i}ker\chi/ker\chi \subseteq Z(G/ker\chi)=Z(\chi)/ker\chi$. Since $|X_{i}|=|Z(\chi)|$, it follows that $Z(\chi)=X_{i}$.

\end{proof}
\end{lemma}

In the next lemma we obtain a necessary and sufficient condition for a character $\chi \in Irr(G)$ to be linear.
\begin{lemma}
\label{lemma4}
Let $\chi \in Irr(G)$. Then $\chi$ is a linear character if and only if $[Z(\chi), G]=G'$.
\begin{proof}
If $\chi$ is a linear character then $[Z(\chi), G]=[G,G]=G'$.\\
Conversely suppose that $[Z(\chi), G]=G'$ where $\chi \in Irr(G)$. Let $\Psi$ be the representation affording the character $\chi$.
Let $g\in G'=[Z(\chi), G]$ so that $g=g_{1}g_{2}\ldots g_{t}$. For each $1\leq i \leq t$ let $g_{i}=x_{i}^{-1}y_{i}^{-1}x_{i}y_{i}$ where $x_{i}\in Z(\chi)$ and $y_{i}\in G$.
Then $\Psi(x_{i})=\alpha_{i}I$ for some $\alpha_{i} \in \mathbb{C}$.
Now $\Psi(g_{i})=\Psi(x_{i}^{-1}y_{i}^{-1}x_{i}y_{i})$\
$=\Psi(x_{i})^{-1}\Psi(y_{i})^{-1}\Psi(x_{i})\Psi(y_{i})$
$=\alpha_{i}^{-1}I\Psi(y_{i})^{-1}\alpha_{i}I\Psi(y_{i})=I$.
Thus $\Psi(g)=I$ and hence 
$ g\in ker\Psi=ker\chi$. This shows that $G'\subseteq ker\chi$. Hence $\chi$ is a linear character.
\end{proof}
\end{lemma}

Since $[Z(\chi),G]<G'$ for any $\chi \in nl(G)$ it makes sense to find the centre of the quotient group $G/[Z(\chi),G]$.

\begin{lemma}
\label{lemma5}
Let $\chi \in Irr(G)$. Then $Z(G/[Z(\chi),G])=Z(\chi)/[Z(\chi), G]$.

\begin{proof}
$Z(\chi)/[Z(\chi), G]\leq Z(G/[Z(\chi),G])$ is obvious.\\
Let $Z/[Z(\chi), G]  =Z(G/[Z(\chi), G])$. Let $\Psi$ be the representation affording the character $\chi$.\\
Let $g[Z(\chi), G]\in Z(G/[Z(\chi), G])$.Then
for any $x\in G$, $gx[Z(\chi), G]=xg[Z(\chi), G]$ and hence $ [g,x]\in [Z(\chi), G]\leq ker\chi=ker\Psi$.
This gives that $\Psi([g,x])=I,  \forall x \in G$ and therefore
$\Psi(g)\Psi(x)=\Psi(x)\Psi(g), \forall x \in G$.
By Lemma 2.25 of \cite{isaacs1994character} we get $\Psi(g)=\alpha I$ for some $\alpha \in \mathbb{C}$.
By Lemma 2.27 of \cite{isaacs1994character} we get that $g\in Z(\chi)$. Hence we conclude that $Z(G/[Z(\chi),G])\leq Z(\chi)/[Z(\chi),G]$. This proves the result. 
\end{proof}
\end{lemma}

In the following lemma we find the centre of a quotient group of certain type of GVZ-group with two character degrees .

\begin{lemma}
\label{lemma6}
Let $G$ be a GVZ group with two character degrees. Let  $N\unlhd G$ be such that $G'\not \subseteq N$ but $[Z(\chi), G]\leq N$ for some $ \chi\in nl(G)$, then $Z(G/N)=Z(\chi)/N$.
\begin{proof}
Let $Z/N=Z(G/N)$. Then $Z=\cap_{\theta \in Irr(G/N)} Z(\theta)$. For every $\theta \in Irr(G/N)$, $N\leq ker \theta$. Then $[Z(\chi),G]\leq N \leq ker \theta$ and hence $Z(\chi)=Z(\theta)$. Thus $Z=\cap_{\theta \in Irr(G/N)} Z(\theta)=Z(\chi)$.\\
\end{proof}
\end{lemma}

Now we state the following definitions.

\begin{defn}
\cite{lewis2019semi}

  A $p$-group $G$ is said to be semi extraspecial if $G$ satisfies the property that for every maximal subgroup $N$ of $Z(G)$, $G/N$ is a extraspecial group.
\end{defn}

   Isoclinicism is an equivalence relation of groups and was first introduced by P. Hall\cite{phall1940} to classify $p$-groups, although it is applicable to all groups.

\begin{defn}
\cite{suzuki1982group}
Two groups $G$ and $H$ are said to be isoclinic if there exist isomorphism $\phi: G/Z(G)\rightarrow H/Z(H)$ such that the mapping $[g_{1},g_{2}]\rightarrow [\phi(g_{1}Z(G)), \phi(g_{2}Z(G))]$ induces an isomorphism between the groups $G'$ and $H'$.
\end{defn}

In \cite{fernandez2001groups} the authors examined the $p$-groups of central type with two character degrees and obtained the following results.

\begin{thm}
\label{thm1}
For a non-abelian $p$-group $G$, the following conditions are equivalent:
  \begin{enumerate}[(i)]

      \item $cd(G)=\{1, |G/Z(G)|^{\frac{1}{2}}\}$.
      \item  $Cl_{G}(x) = xG'$ for all $x \in G \setminus Z(G)$.
      \item $G'=[x,G]$ for all $x \in G \setminus Z(G)$.

      \item  $G$ is isoclinic to a semiextraspecial 
             $p$-group.
      \item $Z(G/N)=Z(G)/N$ for any normal subgroup $N$ of $G$ such that $G'\not \subseteq N$.
  \end{enumerate}
\end{thm}

\section{Main Results}
In the following theorem we determine the structure of some of the factor groups of a GVZ-group with two character degrees.

\begin{thm}
Let $G$ be a  GVZ-group with two character degrees. Then for any $\chi\in nl(G)$, the groups  $G/Z(\chi)$, $G/Z(G)$ and $G'/[Z(\chi),G]$
  are elementary Abelian and $exp\ G/Z(\chi)=exp\ G/Z(G)=exp\ G'/[Z(\chi),G]$.

\begin{proof}
We get that $c(G)=2$. Therefore, $1\neq G'\subseteq Z(G)$. Since the irreducible character $\chi$ is non-linear, we get that $[Z(\chi),G]\not \subseteq G'\subseteq Z(G)$ and hence we can find a prime $p$ and a normal subgroup $N$ of $G$ such that $[Z(\chi),G]\subseteq N \subseteq G'$ and $|G':N|=p$. Since $Z(G/N)=Z(\chi)/N$ and $G/Z(\chi)\cong (G/N)/(Z(\chi)/N)$, we can assume that $N=1$. 
Then $|G'|=p$. Since the nilpotency class of $G$ is $2$, by Lemma 4.4 \cite{isaacs2008finite} we get that exp $G/Z(G)=p$ and hence 
$G/Z(G)$ is elementary Abelian. Since $G/Z(\chi)\cong (G/Z(G))/(Z(\chi)/Z(G))$, it follows that $G/Z(\chi)$ is an elementary Abelian group for any $\chi \in nl(G)$.\\

Let $\chi_{1}, \chi_{2}\in nl(G)$.
Suppose $G/Z(\chi_{1})$ is an elementary Abelian $p_{1}$-group and $G/Z(\chi_{2})$ is an elementary Abelian $p_{2}$-group.
Now for any $x\in G$, $x^{p_{1}}\in Z(\chi_{1}) \subseteq Z(\chi_{1})Z(\chi_{2})$ and this shows that $G/Z(\chi_{1})Z(\chi_{2})$ is an elemantary Abelian $p_{1}$-group. Again $x^{p_{2}}\in Z(\chi_{2})\subseteq Z(\chi_{1})Z(\chi_{2})$ and this shows $G/Z(\chi_{1})Z(\chi_{2})$ is an elemantary Abelian $p_{2}$-group. Hence $p_{1}=p_{2}$. This proves that $exp\ G/Z(\chi)=p$ for all $\chi\in nl(G)$.\\

Hence for any $x\in G$, we get that $x^{p}\in Z(\chi)$ for all $\chi \in nl(G)$. This in turn gives that $x^{p}\in \cap_{\chi \in nl(G)}Z(\chi)=Z(G)$. Thus we get that the group $G/Z(G)$  is elementary Abelian and $exp\ G/Z(G)=p$.\\

Since $G'/[Z(\chi), G]=(G/[Z(\chi), G])'$, it is enough to assume that $[Z(\chi), G]=1$. The group $G'$ is generated by the elements of the form $[x,y]$, $x,y\in G$. Since the nilpotency class of $G$ is $2$ and $x^{p}\in Z(G)$ for all $g\in G$, we get that $[x,y]^{p}=[x^{p},y]=1$. Thus exponent of every element of $G'$ is $p$ and consequently $G'$ is an elementary Abelian group.\\

\end{proof}
\end{thm}

From Lemma \ref{lemma4} we get that $\chi\in Irr(G)$ is a non-linear character if and only if $[Z(\chi), G]<G'$. In the following theorem we get a similar result for the  irreducible characters of the factor groups of a GVZ-group with two character degrees.

\begin{thm}
Let $G$ be a GVZ group with two character degrees. Then  $\chi \in nl(G/N)$  if and only if $N[Z(\chi),G]<NG'$.

\begin{proof}
First we assume that there is a $\chi \in nl(G/N)$ . Suppose $N[Z(\chi),G]=NG'$. By Lemma \ref{lemma3} we get that $[Z(\chi), G]\leq ker\chi$. Again $N\leq ker\chi$.Thus $G'\leq NG'=N[Z(\chi),G]\leq ker\chi$. This gives that  $\chi$ is a linear character. This contradiction proves that $N[Z(\chi),G]<NG'$.\\
Conversely assume that $N[Z(\chi),G]<NG'$ for $\chi \in Irr(G/N)$. We choose a non-principal irreducible character $\lambda \in Irr(NG'/N[Z(\chi),G])$. If possible assume that all the irreducible constituent characters of $\lambda^{G}$ are linear. Then $G'\subseteq ker(\lambda^{G}) \subseteq ker\lambda$ and hence $NG'\subseteq ker\lambda$. This gives that $\lambda$ is the principal character of the group $NG'/N[Z(\chi),G]$. Thus there exists a non-linear irreducible character $\phi $ of $G$ which is a constituent of $\lambda^{G}$. Then $N[Z(\chi),G]\leq ker\lambda \leq ker\phi$. The  fact $[Z(\chi),G] \leq ker\phi$ gives that 
$Z(\phi)=Z(\chi)$. Since $\phi$ is a non-linear irreducible character we must get that $\chi$ is a non-linear irreducible character of $G/N$.
\end{proof}
\end{thm}

In the following theorem we characterize the GVZ-groups with two character degrees by the orders of the centres of some factor groups of the group.

\begin{thm}
Let $G$ be a  a GVZ group. Then the following conditions are equivalent:
  \begin{enumerate}[(i)]
   \item $|cd(G)|=2$.

   \item There is a collection $\{X_{1}, X_{2}, \ldots , X_{n}\}$ of normal subgroups of $G$ of equal order such that  every normal subgroup $N$ of $G$ with $G'\not \subseteq N$ and $[Z(\chi),G ]\subseteq N$ for some $\chi \in \ nl(G)$ satisfies $Z(G/N)=X_{i}/N$ for some $1\leq i \leq n$.

   \item For any two normal subgroups $N$ and $M$ of $G$ such that $G'\not \subseteq N,M$ and $[Z(\chi), G] \subseteq N$, $[Z(\phi), G] \subseteq M$ for some $\chi, \phi \in \ nl(G)$ if we define $Z_{N}/N=Z(G/N)$, $Z_{M}/M=Z(G/M)$ then $|Z_{N}|=|Z_{M}|$.
   
Moreover the collection $\{X_{1}, X_{2}, \ldots , X_{n}\}$ of normal subgroups in (ii) is the collection of centres of non-linear irreducible characters of $G$.
  \end{enumerate}

 \begin{proof}

$(i)\Rightarrow (ii)$
Since $G$ is a GVZ-group with $|cd(G)|=2$, we get that $|Z(\chi)|$ is constant for all $\chi \in nl(G)$. Let $\{X_{1}, X_{2}, \ldots , X_{n}\}$ be the collection of centres of non-linear irreducible characters of $G$. Let $N$ be a normal subgroup of $G$ such that $G'\not \subseteq N$ and $[Z(\chi), G] \subseteq N$ for some $\chi \in \ nl(G)$. Let $Z(\chi)=X_{i}$  . Then by Lemma \ref{lemma6} $Z(G/N)=Z(\chi)/N=X_{i}/N$. \\

$(ii)\Rightarrow (iii)$
Let $N$ and $M$ be normal subgroups of $G$ such that $G'\not \subseteq N,M$ and $[Z(\chi), G] \subseteq N$, $[Z(\phi), G] \subseteq M$ for some $\chi, \phi \in \ nl(G)$. If $Z_{N}/N=Z(G/N)$, $Z_{M}/M=Z(G/M)$ we get that $Z_{N}=X_{i}$ and $Z_{M}=X_{j}$ for some $i,j$. Since $|X_{i}|=|X_{j}|$ we get that $|Z_{N}|=|Z_{M}|$.\\

$(iii)\Rightarrow (i)$
  Let $\chi,\phi \in \ nl(G)$.
  Then $G'\not \subseteq ker\chi, ker \phi$ and $[Z(\chi),G]\leq ker\chi$, $[Z(\phi),G]\leq ker\phi$. Since $Z(\chi)/ker\chi=Z(G/ker\chi)$ and $Z(\phi)/ker\phi=Z(G/ker\phi)$, it follows that $|Z(\chi)|=|Z(\phi)|$. Since $G$ is a GVZ-group it follows that $\chi(1)^{2}=|\frac{G}{Z(\chi)}|$ is constant and hence we get that $|cd(G)|=2$.\\

  For every $\chi \in nl(G)$, we get that $G'\not \subseteq ker\chi$ and $[Z(\chi),G]\subseteq ker \chi$. Thus by (ii) we get that $Z(\chi)/ker \chi=Z(G/ker \chi)=X_{i}/ker\chi$ for some $1\leq i \leq n$ and hence $Z(\chi)=X_{i}$. This shows that $\{X_{1}, X_{2}, \ldots , X_{n}\}$ is the collection of the centres of the non-linear irreducible characters of $G$.

      \end{proof}
\end{thm}

In our last theorem we obtain some characterizations of non-Abelian $p$-groups which are GVZ-groups with two character degrees. The results of this theorem generalize the results of Theorem \ref{thm1}.

\begin{thm}
Let $G$ be a non-Abelian $p$ group with $|cd(G)|=2$ . Then the following conditions are equivalent:
\begin{enumerate}[(i)]
    \item  $G$ is a GVZ group.
   
    %\item $Cl(x[Z(\chi), G])=xG'/[Z(\chi), G]$ in the group $G/[Z(\chi),G]$ for all $x \not \in  Z(\chi)$, $\chi \in \ nl(G)$.

     \item $G/[Z(\chi),G]$ is isoclinic to a semi extraspecial $p$-group for any $\chi \in nl(G)$.
     
    \item $Z(G/N)=Z(\chi)/N$ for any normal subgroup $N$ of $G$ such that $G'\not \leq N$ and $[Z(\chi),G]\leq N$ for some $\chi \in \ nl(G)$.
    
\end{enumerate}

\begin{proof}

$(i)\Rightarrow (ii)$ $G$ is a GVZ-group with two character degrees.  Let
$\overline{G}=G/[Z(\chi), G]$. By Lemma \ref{lemma5} we get that $\overline{G}$ is a GVZ-group with $cd(\overline{G})=\{1, |\frac{\overline{G}}{Z(\overline{G})}|^{\frac{1}{2}}\}$. Hence by Theorem \ref{thm1} we get that $G/[Z(\chi), G]$ is isoclinic to a semi extraspecial $p$-group.
\\

$(ii)\Rightarrow (i)$  Let $\chi \in nl(G)$ and 
$\overline{G}=G/[Z(\chi), G]$. Since $\overline{G}$ is isoclinic to a semi-extraspecial $p$-group, it follows from Theorem \ref{thm1}   that $cd(\overline{G})=\{1, |\overline{G}:Z(\overline{G})|^{\frac{1}{2}}\}$.Since $Z(G/[Z(\chi),G])=Z(\chi)/[Z(\chi),G]$, it follows that $cd(\overline{G})=\{1,|G:Z(\chi)|^{\frac{1}{2}}\}$. Since $|cd(G)|=2$ we get that $cd(G)=\{1,|G:Z(\chi)|^{\frac{1}{2}}\}$. By Corollary 2.30 of  \cite{isaacs1994character} it follows that $G$ is a GVZ group.
\\

$(i)\Rightarrow (iii)$ 
Let $N$ be a normal subgroup of $G$ such that $G'\not \leq N$ and $[Z(\chi),G]\leq N$ where $\chi \in \ nl(G)$. 
 Since $G$ is a GVZ-group with two character degrees, we get that $N\leq Z(\chi)$.
 Let $\overline{G}=G/[Z(\chi), G]$ and $\overline{N}=N/[Z(\chi), G]$. As shown in the proof of $(i)\Rightarrow (ii)$ we get that $cd(\overline{G})=\{1, |\overline{G}:Z(\overline{G})|^{\frac{1}{2}}\}$. Then by Theorem \ref{thm1} we get that $Z(G/N)\cong Z(\overline{G}/\overline{N})=Z(\overline{G})/\overline{N} \cong Z(\chi)/N$.
 Since $[Z(\chi),G]\leq N$ it follows that $Z(\chi)/N \subseteq Z(G/N)$ and hence $Z(G/N)=Z(\chi)/N$.
\\

 $(iii)\Rightarrow (i)$ Suppose that $Z(G/N)=Z(\chi)/N$ for any normal subgroup $N$ of $G$ such that $G'\not \leq N$ and $[Z(\chi),G]\leq N$ where $\chi \in \ nl(G)$. For any $\chi \in nl(G)$ we consider the normal subgroup $N=ker\chi$ of $G$. Then $G'\not \subseteq N$, $[Z(\chi),G]\leq N$ and $Z(G/N)=Z(\chi)/N$. 
 If  $\overline{G}=G/[Z(\chi), G]$
 by  Theorem \ref{thm1} we get that $cd(\overline{G})=\{1, |\overline{G}: Z(\overline{G})|^{\frac{1}{2}}\}$. Since $|\overline{G}: Z(\overline{G})|=|G:Z(\chi)|$ and $|cd(G)|=2$, we get that $cd(G)=\{1, |G:Z(\chi)|^{\frac{1}{2}}\}$. Thus $\chi(1)^{2}=|G:Z(\chi)|$. Hence $G$ is a GVZ-group.
\end{proof}

\end{thm}

\end{document}